\documentclass[12pt]{amsart}
\usepackage{enumerate}
\usepackage{amssymb}
\usepackage{amsmath}
\usepackage{amsthm}
\newcommand{\snug}{\unskip\kern-\mathsurround}

\newcommand{\Z}{{\mathbb Z}}
\newcommand{\F}{{\mathbb F}}

\newtheorem{theorem}{Theorem}
\newtheorem{lemma}[theorem]{Lemma}

\theoremstyle{definition}

\begin{document}
\title{Residually Torsion-Free Nilpotent One Relator Groups}
\author{John Labute}
\address{Department of Mathematics and Statistics, McGill University, Burnside Hall, 805 Sherbrooke Street West, Montreal QC H3A 0B9, Canada}
\email{labute@math.mcgill.ca}
\thanks{I would like to thank Roman Mikhailov for pointing out the work of Azarov.}
\begin{abstract} We show that the group $G=\langle x_1,..x_m,y_1,..y_n \mid u=v\rangle$ is residually torsion-free nilpotent if $v\in\langle y_1,..y_n\rangle$, $v\ne1$, $u\in A= \langle x_1,..x_m\rangle$, $u\in \gamma_d(A)$, $u$ not a proper power mod $\gamma_{d+1}(A)$, where $\gamma_k(A)$ is the $k$-th term of the lower central series of $A$.
\end{abstract}
\date{March 8, 2015}
\maketitle
\textbf{}\hfill{\it In memory of Gilbert Baumslag}
\bigskip

In \cite{Az} Azarov proves that the group
$$
G=\langle x_1,..x_m,y_1,..y_n \mid u=v\rangle
$$
is residually a finite $p$-group for any prime $p$ if $u\in A=\langle x_1,..x_m\rangle$ is not a proper power and $v\in B=\langle y_1,..y_n\rangle$, $v\ne1$.

Let $\gamma_k(A)$ be the $k$-th term of the lower central series of $A$. If we strengthen the condition on $u$ by requiring that $u\in \gamma_d(A)$ but not a proper power mod $\gamma_{d+1}(A)$ we obtain the following result which extends a result of Baumslag and Mikhailov (cf. \cite{BM}, Theorem 5).
\begin{theorem}
If $u\in \gamma_d(A)$ and $u$ is not a proper power mod $\gamma_{d+1}(A)$ there exists a central series $(G_i)$ in $G$ such that the quotients $G_i/G_{i+1}$ are torsion free and the intersection of the groups $G_i$ is $1$. In particular $G$ is residually torsion-free nilpotent.
\end{theorem}
A sequence of subgroups $G_i(i\ge1)$ is said to be a central series for $G$ if $G_1=G$, $G_{i+1}\subseteq G_i$, $[G_i,G_j]\subseteq G_{i+j}$. To construct the required central series we will make use of the Magnus embedding of the free group $F=A*B$ into the Magnus algebra $M$ of formal power series in the non-commuting variables $X_1,..X_m,Y_1,..Y_n$ with integer coefficients which sends $x_i$ into $1+X_i$ and $y_i$ into $1+Y_i$. 

Let $e$ be an integer greater than the integer $d$ in Theorem 1. If
$$
f=\sum a_{i_1\ldots i_k}Z_{i_1}\cdots Z_{i_k}\in M
$$
is a sum of distinct monomials $Z_{i_1}\cdots Z_{i_k}$ with $a_{i_1\ldots i_k}\in \Z$ $(k\ge0)$ and $Z_i$ an element of $\{X_1,\ldots,X_m,Y_1,\ldots,Y_n\}$, we define a valuation $v$ on $M$ by
$$
v(f)= {\rm min} \{a + eb\mid a_{i_1\ldots i_k}\ne0\}
$$
where $a$ and $b$ are respectively the number of the $X_i$ and $Y_i$ in $Z_{i_1}\cdots Z_{i_k}$. In particular we have $v(X_i)=1, v(Y_i)=e$. By convention $v(0)=\infty$.

For $i\ge1$ let $M_i=\{f\in M\mid v(f)\ge i$\} and let $F_i=F\cap(1+M_i)$ where we have identified $F$ with its image in $M$. Then $(F_i)$ is a central series for $F$ (cf. \cite{La3}, sect.~ 2). (Note that if $e=1$ in the definition of $v$ we have $F_i=\gamma_i(F)$ by a deep result of Magnus, cf. \cite{La3}, section 3.) If $G_i$ is the image of $F_i$ in $G$ we will show that $(G_i)$ is the required central series.
\medskip

\begin{lemma}
For $i\ge e$ we have $F_i\subseteq \gamma_{[i/e]}(F)$ which shows that $G_i\subseteq \gamma_{[i/e]}(G)$ and hence that $\bigcap G_i=\bigcap\gamma_i(G)=1$.
\end{lemma}
This follows from the fact that if $i=a+eb$ then 
$$
a+b=(ea+eb)/b\ge (a+eb)/e=i/e
$$ 
and the fact that $G$ is residually nilpotent by Azarov's Theorem.
\medskip

In order to show that the quotients $G_i/G_{i+1}$ are torsion-free we have to bring into play the Lie ring structure on the graded ring $L(G)=\oplus_i (G_i/G_{i+1})$. The Lie ring $L(F)=\oplus_i(F_i/F_{i+1})$ is free on $\xi_1,..\xi_m,\eta_1,..\eta_n$ where $\xi_i$ is the image of $X_i$ in $F_1/F_2$ and $\eta_j$ is the image of $Y_j$ in $F_e/F_{e+1}$ (cf. \cite{La3}, section 3). The problem is to determine the kernel of the canonical surjection $L(F)\rightarrow L(G)$. Let $\rho$ be the image of $r=uv^{-1}$ in $F_d/F_{d+1}$. Since $v\in F_{d+1}$, we have that $\rho$ is also the image of $u$ in $F_d/F_{d+1}$.

\begin{theorem}
The Lie ring $\mathfrak g=L(F)/(\rho)$ is torsion-free.
\end{theorem}

\begin{theorem}
The kernel of $L(F)\rightarrow L(G)$ is $\mathfrak r =(\rho)$.
\end{theorem}

To prove Theorem 3 we note that since $\rho$ is not a proper multiple of an element of $L(F)$ the Lie algebra $(L(F)/(\rho))\otimes\F_p$ is a graded Lie algebra over the finite field $\F_p$ defined by a single non-zero relator of degree $d$. In \cite{La1}, th\'eor\`eme 2 we prove that the homogeneous component of degree $n$ of this graded algebra has a finite dimension which depends only on $n$ and $d$. Since $p$ is an arbitrary prime this proves that the homogeneous components of $\mathfrak g=L(F)/(\rho)$ are torsion free. The Theorem of Birkhoff-Witt then shows that the enveloping algebra $U$ of $\mathfrak g$ has no zero-divisors which can be used to prove that, via the adjoint representation, $\mathfrak r/[\mathfrak r,\mathfrak r]$ is a free $U$ module generated by the image of $\rho$ (cf. \cite{La1}, th\'eor\`eme 1). This fact is the key to proving Theorem 4. The proof given in \cite{La2} can be easily be adapted to the valuation $v$ used here. For details cf. \cite{La3}, sections 2 and 3.

\end{document}